\documentclass[10pt,twoside,reqno]{amsart}
\usepackage{amssymb}
\usepackage{multirow}
\calclayout

\numberwithin{equation}{section}
\makeatletter

\renewcommand{\@secnumfont}{\bfseries}

\renewcommand{\section}{\@startsection{section}{1}%
  {0mm}{.7\linespacing\@plus\linespacing}{.5\linespacing}
  {\normalfont\bfseries\centering}}

\newcommand{\bibsection}{\@startsection{section}{1}%
  {0mm}{.7\linespacing\@plus\linespacing}{.5\linespacing}
  {\normalfont\scshape\centering}}

\renewcommand{\@biblabel}[1]{#1.}

\newtheorem{thm}{\bf Theorem}[section]

\begin{document}

\vspace{1.3cm}

\title {A note on degenerate Stirling numbers of the first kind}

\author{Taekyun Kim}
\address{Department of Mathematics, Kwangwoon University, Seoul 139-701, Republic of Korea}
\email{tkkim@kw.ac.kr}

\author{Dae San Kim}
\address{Department of Mathematics, Sogang University, Seoul 121-742, Republic of Korea}
\email{dskim@sogang.ac.kr}

\begin{abstract} 
Recently, the degenerate Stirling numbers of the first kind were introduced. In this paper, we give some formulas for the degenerate Stirling numbers of the first kind in the terms of the complete Bell polynomials with higher-order harmonic number arguments. Also, we derive an identity connecting the degenerate Stirling numbers of the first kind and the degenerate derangement numbers by using probabilistic method.
\end{abstract}
\maketitle
\bigskip
\medskip

\section{Introduction}

The numbers $S_1(n,k)$ and $S_2(n,k)$ are, in the notation of Riordan [13,14], Stirling numbers of the first kind and of the second kind, respectively and they are given by:
\begin{equation}\begin{split}\label{01}
(x)_n = \sum_{k=0}^n S_1(n,k) x^k,\,\,(n \geq 0),
\end{split}\end{equation}
\begin{equation}\begin{split}\label{02}
x^n = \sum_{k=0}^n S_2(n,k) (x)_k ,\,\,(n \geq 0),\quad (\text{see}\;[1,2,3]).
\end{split}\end{equation}

\noindent where $(x)_n=x(x-1)\cdots(x-(n-1))$, $(n \geq 1)$, and $(x)_0=1$. 

The unsigned Stirling numbers of the first kind are defined as
\begin{equation}\begin{split}\label{03}
<x>_n = \sum_{k=0}^n S(n,k) x^k,\,\,(n \geq 0),\quad (\text{see}\;[13]),
\end{split}\end{equation}
where $<x>_0=1$, $<x>_n = x(x+1)\cdots(x+n-1), (n \geq 1)$.

From \eqref{01} and \eqref{03}, we note that
\begin{equation}\begin{split}\label{04}
S(n,k) = (-1)^{n-k} S_1(n,k),\,\,(n \geq k \geq 0).
\end{split}\end{equation}

Let $\lambda$ be a real number. Then the degenerate Euler polynomials are defined by the generating function 
\begin{equation}\begin{split}\label{05}
\frac{2}{(1+\lambda t)^{\frac{1}{\lambda }}+1} (1+\lambda t)^{\frac{x}{\lambda }} = \sum_{n=0}^\infty       \mathcal{E}_{n,\lambda }(x) \frac{t^n}{n!},\quad (\text{see}\;[1]).
\end{split}\end{equation}

When $x=0$, $\mathcal{E}_{n,\lambda }=\mathcal{E}_{n,\lambda }(0)$ are called the degenerate Euler numbers.

From \eqref{05}, we note that $\lim_{\lambda \rightarrow 0} \mathcal{E}_{n,\lambda }(x) = E_n(x)$, $(n \geq 0)$, where $E_n(x)$ are the Euler polynomials given by 
\begin{equation*}\begin{split}
\frac{2}{e^t+1} e^{xt} = \sum_{n=0}^\infty     E_n(x)   \frac{t^n}{n!},\quad (\text{see}\;[2,12,13,14]).
\end{split}\end{equation*}

By \eqref{01} and \eqref{03}, we easily get
\begin{equation}\begin{split}\label{06}
\frac{1}{k!} \Big( \log(1+t) \Big)^k = \sum_{n=k}^\infty S_1(n,k) \frac{t^n}{n!},\quad (\text{see}\;[13]),
\end{split}\end{equation}
and
\begin{equation}\begin{split}\label{07}
\frac{1}{k!} \Big( \log \big( \frac{1}{1-t} \big) \Big)^k = \sum_{n=k}^\infty S(n,k) \frac{t^n}{n!}.
\end{split}\end{equation}

For $\lambda  \in \mathbb{R}$, the degenerate falling factorial sequence is defined as
\begin{equation}\begin{split}\label{08}
(x)_{0,\lambda }=1, \,\,(x)_{n,\lambda }= x(x-\lambda )(x-2\lambda ) \cdots (x-(n-1)\lambda ),\,\,(n \geq 1).
\end{split}\end{equation}

Note that $\lim_{\lambda \rightarrow 1} (x)_{n,\lambda } = (x)_n$, $\lim_{\lambda  \rightarrow 0 } (x)_{n,\lambda } = x^n$. In [5], the degenerate Stirling numbers of the first kind were defined by Kim as
\begin{equation}\begin{split}\label{09}
(x)_{n,\lambda } = \sum_{k=0}^n S_{1,\lambda }(n,k) x^k,\,\,(n \geq 0).
\end{split}\end{equation}

K\"{o}lbig gave the following formula:
\begin{equation}\begin{split}\label{10}
\frac{d^r}{dx^r} e^{f(x)} = e^{f(x)} Bel_r \left( f^{(1)}(x), f^{(2)}(x), \cdots, f^{(r)}(x) \right),\quad (\text{see}\;[10,11]),
\end{split}\end{equation}

\noindent where $f^{(r)}(x) = \left( \frac{d}{dr} \right)^r f(x)$ and $Bel_r(x_1,x_2,\cdots,x_r)$ are the complete Bell polynomials given by
\begin{equation*}\begin{split}
Bel_r (x_1,x_2,\cdots,x_r) = \sum_{k_1+2k_2+\cdots+rk_r=r} {r \choose k_1,\cdots,k_r} \left( \frac{x_1}{1!} \right)^{k_1} \left( \frac{x_2}{2!} \right)^{k_2} \cdots  \left( \frac{x_r}{r!} \right)^{k_r}.
\end{split}\end{equation*}

The complete Bell polynomials are also given by the exponential generating function 
\begin{equation}\begin{split}\label{11}
\exp \left( \sum_{j=1}^\infty   x_j     \frac{t^j}{j!} \right) = \sum_{n=0}^\infty  Bel_n(x_1,x_2,\cdots,x_n)      \frac{t^n}{n!},\quad (\text{see}\;[4]).
\end{split}\end{equation}

In [4], the complete degenerate Bell polynomials are defined by 

\begin{equation}\begin{split}\label{12}
\exp \left( \sum_{j=1}^{\infty}	 x_j (1)_{j,\lambda } \frac{t^j}{j!} \right) =\sum_{n=0}^\infty       Bel_n^{(\lambda )}(x_1,x_2,\cdots,x_n) \frac{t^n}{n!},
\end{split}\end{equation}
where
\begin{equation*}\begin{split}
Bel_n^{(\lambda )}(x_1,\cdots,x_n) = \sum_{k_1+2k_2\cdots+nk_n=n} {n \choose k_1,\cdots,k_n} \left( \frac{x_1 (1)_{1,\lambda }}{1!} \right)^{k_1} \cdots \left( \frac{x_n (1)_{n,\lambda }}{n!} \right)^{k_n}.
\end{split}\end{equation*}

From \eqref{12}, we note that
\begin{equation*}\begin{split}
Bel_n^{(\lambda )} (1,1,\cdots,1) = Bel_{n,\lambda }
\end{split}\end{equation*}
are the degenerate Bell numbers given by
\begin{equation}\begin{split}\label{13}
e^{(1+\lambda t)^{\frac{1}{\lambda }}-1} = \sum_{n=0}^\infty    Bel_{n,\lambda }    \frac{t^n}{n!},\quad (\text{see}\;[9]).
\end{split}\end{equation}

Here the finite Hurwitz-type functions of order $r$ are defined as
\begin{equation}\begin{split}\label{14}
\sum_{k=0}^{n-1} \frac{1}{(k+x)^r} = H_n^{(r)}(x),\,\,(n, r \in \mathbb{N}).
\end{split}\end{equation}
When $x=1$, $H_n^{(r)}=H_n^{(r)}(1) = \sum_{k=0}^{n-1} \frac{1}{(k+1)^r} = \sum_{k=1}^n	\frac{1}{k^r}$ are the harmonic numbers of order $r$.

From \eqref{09}, we can derive the unsigned Stirling numbers of the first kind as follows:
\begin{equation}\begin{split}\label{15}
<x>_{n,\lambda } = \sum_{k=0}^n S_{\lambda }(n,k) x^k,\quad (\text{see}\;[5,6]),
\end{split}\end{equation}
where $<x>_{0,\lambda }=1$, $<x>_{n,\lambda } = x(x+\lambda ) \cdots (x+(n-1)\lambda ), (n \geq 1)$. By \eqref{09} and \eqref{15}, we get
\begin{equation}\begin{split}\label{16}
S_\lambda  (n,k) = (-1)^{n-k} S_{1,\lambda }(n,k),\,\,(n \geq k \geq 0).
\end{split}\end{equation}

In this paper, we give some formulas for the degenerate Stirling numbers of the first kind in the terms of the complete Bell polynomials with higher-order harmonic number arguments. Also, we derive an identity connecting the degenerate Stirling numbers of the first kind and the degenerate derangement numbers by using probabilistic method.

\section{Degenerate Stirling numbers of the first kind}

From \eqref{01}, we have
\begin{equation}\begin{split}\label{16-1}
(1+\lambda t)^{\frac{x}{\lambda }} &= \sum_{k=0}^\infty (x)_{k,\lambda } \frac{t^k}{k!} = \sum_{k=0}^\infty \sum_{n=0}^k S_{1,\lambda }(k,n)x^n \frac{t^k}{k!} \\
&= \sum_{n=0}^\infty \left ( \sum_{k=n}^\infty  S_{1,\lambda }(k,n) \frac{t^k}{k!} \right) x^n.
\end{split}\end{equation}

On the other hand
\begin{equation}\begin{split}\label{17}
(1+\lambda t)^{\frac{x}{\lambda }} &= \sum_{n=0}^\infty \left( \frac{x}{\lambda } \right)^n \frac{1}{n! } \Big( \log(1+\lambda t) \Big)^n \\
&= \sum_{n=0}^\infty \frac{1}{n!} \Big( \log(1+\lambda t)^{\frac{1}{\lambda }} \Big)^n x^n.
\end{split}\end{equation}

By \eqref{16-1} and \eqref{17}, we get
\begin{equation}\begin{split}\label{18}
\frac{1}{n!} \Big( \log(1+\lambda t)^{\frac{1}{\lambda }} \Big)^n = \sum_{k=n}^\infty  S_{1,\lambda }(k,n) \frac{t^k}{k!}.
\end{split}\end{equation}

It is easy to show that
\begin{equation}\begin{split}\label{19}
(-x)_{n,\lambda } = (-1)^n <x>_{n,\lambda },\,\,(n \geq 1).
\end{split}\end{equation}

Thus, by \eqref{19}, we see again that
\begin{equation}\begin{split}\label{20}
S_\lambda  (n,k) = (-1)^{n-k} S_{1,\lambda }(n,k),\,\,(n \geq k \geq 0).
\end{split}\end{equation}

From \eqref{19}, we also note that
\begin{equation}\begin{split}\label{21}
\frac{1}{k!} \left( \log \Big( \frac{1}{(1-\lambda t)^{\frac{1}{\lambda }}} \Big) \right)^k = \sum_{n=k}^\infty S_\lambda (n,k) \frac{t^n}{n!}.
\end{split}\end{equation}

By \eqref{05}, we easily get
\begin{equation}\begin{split}\label{22}
2 \sum_{l=0}^{m-1} (1+\lambda t)^{\frac{l}{\lambda }} (-1)^l = \sum_{n=0}^\infty \Big( \mathcal{E}_{n,\lambda } + \mathcal{E}_{n,\lambda }(m) \Big) \frac{t^n}{n!},
\end{split}\end{equation}
where $m \in \mathbb{N}$ with $m \equiv 1$ (mod 2).

On the other hand
\begin{equation}\begin{split}\label{23}
2\sum_{l=0}^{m-1}(1+\lambda t)^{\frac{l}{\lambda }}(-1)^l &= \sum_{n=0}^\infty \left( 2 \sum_{l=0}^{m-1} (l)_{n,\lambda }(-1)^l \right) \frac{t^n}{n!} \\
&= \sum_{n=0}^\infty \left( 2 \sum_{k=0}^n \sum_{l=0}^{m-1} (-1)^l S_{1,\lambda }(n,k) l^k     \right)  \frac{t^n}{n!}.
\end{split}\end{equation}

Thus, by \eqref{22} and \eqref{23}, we get
\begin{equation}\begin{split}\label{24}
\mathcal{E}_{n,\lambda } + \mathcal{E}_{n,\lambda }(m) = 2 \sum_{k=0}^n \sum_{l=0}^{m-1} (-1)^l S_{1,\lambda }(n,k) l^k,
\end{split}\end{equation}
where $n \geq 0$, and $m \in \mathbb{N}$, with $m \equiv 1$ (mod 2).

We summarize this as a theorem.

\begin{thm}
Let  $S_{1,\lambda }(n,k)$ be the degenerate Stirling numbers of the first kind given in \eqref{09} or \eqref{18}. Then we have
\begin{equation}\begin{split}
\mathcal{E}_{n,\lambda } + \mathcal{E}_{n,\lambda }(m) = 2 \sum_{k=0}^n \sum_{l=0}^{m-1} (-1)^l S_{1,\lambda }(n,k) l^k,\nonumber
\end{split}\end{equation}
where $n \geq 0$, and $m \in \mathbb{N}$, with $m \equiv 1$ (mod 2).
\end{thm}

Now, we observe that
\begin{equation}\begin{split}\label{25}
\frac{d}{dx} \frac{<x>_{n,\lambda }}{x} &= \frac{d}{dx} e^{\big(\log(<x>_{n,\lambda }-\log x \big)}\\
&= \frac{<x>_{n,\lambda }}{x} \frac{d}{dx} \left( \sum_{k=0}^{n-1} \log(x+k\lambda ) - \log x \right) \\
&= \frac{<x>_{n,\lambda }}{x} \sum_{k=1}^{n-1} \frac{1}{x+k\lambda }.
\end{split}\end{equation}

Let $T_\lambda (x) = \sum_{k=1}^{n-1} \frac{1}{x+k\lambda }.$ 
Then we have
\begin{equation}\begin{split}\label{26}
\frac{d^r}{dx^r} T_\lambda (x) = T_\lambda ^{(r)}(x) = (-1)^r r! \sum_{k=1}^{n-1} \frac{1}{(x+k\lambda )^{r+1}}.
\end{split}\end{equation}

By \eqref{26}, we get
\begin{equation}\begin{split}\label{27}
T_\lambda ^{(r)}(0) = (-1)^r r! \lambda ^{-r-1} \sum_{k=1}^{n-1} \frac{1}{k^{r+1}} = (-1)^r r! \lambda ^{-r-1} H_{n-1}^{(r+1)}.
\end{split}\end{equation}

From \eqref{10}, we have
\begin{equation}\begin{split}\label{28}
\frac{d^r}{dx^r} \frac{<x>_{n,\lambda }}{x} = \frac{<x>_{n,\lambda }}{x} Bel_r \left(T_\lambda (x), T_\lambda ^{(1)}(x), \cdots, T_\lambda ^{(r-1)}(x) \right).
\end{split}\end{equation}

By \eqref{28}, we get
\begin{equation}\begin{split}\label{29}
& \frac{d^r}{dx^r} \frac{<x>_{n,\lambda }}{x} \Bigg|_{x=0}\\
& = (n-1)!\lambda^{n-1} Bel_r \left( 
\frac{1}{\lambda }H_{n-1}^{(1)}, \frac{(-1)1!}{\lambda^2}H_{n-1}^{(2)},\cdots,\frac{(-1)^{r-1}(r-1)!}{\lambda ^r}H_{n-1}^{(r)}
\right).
\end{split}\end{equation}

The equation \eqref{29} is equivalent to
\begin{equation}\begin{split}\label{30}
&\frac{d^r}{dx^r} \frac{<x>_{n,\lambda }}{x}\Bigg|_{x=0}\\
& = (n-1)! \lambda ^{n-1-r} Bel_r \left( H_{n-1}^{(1)}, (-1)1!H_{n-1}^{(2)}, \cdots, (-1)^{r-1}(r-1)! H_{n-1}^{(r)} \right).
\end{split}\end{equation}

On the other hand,
\begin{equation}\begin{split}\label{31}
\frac{d^r}{dx^r} \frac{<x>_{n,\lambda }}{x} &= \frac{d^r}{dx^r} \sum_{l=0}^n S_\lambda (n,l) x^{l-1}\\
&= \sum_{l=0}^n S_\lambda (n,l) (l-1) (l-2) \cdots (l-r)x^{l-1-r}.
\end{split}\end{equation}

Thus, we have
\begin{equation}\begin{split}\label{32}
\frac{d^r}{dx^r} \frac{<x>_{n,\lambda }}{x} \Bigg|_{x=0} = S_\lambda (n,r+1)r!,
\end{split}\end{equation}
where $n \geq r+1 \geq 1$. 

Therefore, by \eqref{30} and \eqref{32}, we obtain the following equation.
\begin{equation}\begin{split}\label{33}
&S_\lambda (n,r+1) \\
&= \frac{(n-1)!}{r!} \lambda ^{n-1-r} Bel_r \left( H_{n-1}^{(1)},(-1)1! H_{n-1}^{(2)},\cdots, (-1)^{r-1}(r-1)! H_{n-1}^{(r)} \right).
\end{split}\end{equation}

In particular, by replacing $n$ by $n+1$, we get
\begin{equation}\begin{split}\label{34}
S_\lambda (n+1,r+1) = \frac{n!}{r!} \lambda ^{n-r} Bel_r \left( H_n^{(1)}, (-1)1! H_n^{(2)},\cdots (-1)^{r-1}(r-1)! H_n^{(r)} \right),
\end{split}\end{equation}
where $n \geq r \geq 0$.

Thus we obtain the following theorem.

\begin{thm}
Let $S_\lambda (n,r)$ be the unsigned degenerate Stirling numbers of the first kind given in \eqref{15} or \eqref{21}. Then, for $n \geq r \geq 0$, we have
\begin{equation}\begin{split}
S_\lambda (n+1,r+1) = \frac{n!}{r!} \lambda ^{n-r} Bel_r \left( H_n^{(1)}, (-1)1! H_n^{(2)},\cdots (-1)^{r-1}(r-1)! H_n^{(r)} \right),\nonumber
\end{split}\end{equation}
where $H_n^{(i)}= \sum_{k=1}^{n}\frac{1}{k^i}$ are the harmonic numbers of order $i$.
\end{thm}

The exponential partial Bell polynomials are the ones given by 
\begin{equation}\begin{split}\label{35}
\frac{1}{k!} \left( \sum_{m=1}^\infty x_m \frac{t^m}{m!} \right)^k = \sum_{n=k}^\infty B_{n,k} (x_1,x_2,\cdots,x_{n-k+1}) \frac{t^n}{n!},\,\,(k \geq 0),
\end{split}\end{equation}
where 
\begin{equation*}\begin{split}
&B_{n,k} (x_1,x_2,\cdots,x_{n-k+1})\\
&= \sum_{\substack{i_1+\cdots+i_{n-k+1}=k\\i_1+2i_2+\cdots+(n-k+1)i_{n-k+1}=n}} {n \choose i_1,\cdots,i_{n-k+1}} \left( \frac{x_1}{1!} \right)^{i_1}  \cdots \left( \frac{x_{n-k+1}}{(n-k+1)!} \right)^{i_{n-k+1}}.
\end{split}\end{equation*}

Recently, Kim defined the degenerate Stirling numbers of the second kind as follows:
\begin{equation}\begin{split}\label{36}
\frac{1}{k!} \Big( (1+\lambda t)^{\frac{1}{\lambda }}-1 \Big)^k = \sum_{n=k}^\infty S_{2,\lambda }(n,k) \frac{t^n}{n!},\quad (\text{see}\;[6]).
\end{split}\end{equation}

From \eqref{36}, the exponential partial $\lambda $-Bell polynomials are considered by Kim as
\begin{equation}\begin{split}\label{37}
\frac{1}{k!} \left( \sum_{i=1}^\infty (1)_{i,\lambda }x_i \frac{t^i}{i!}\right)^k = \sum_{n=k}^\infty B_{n,k}^{(\lambda )} (x_1,\cdots,x_{n-k+1}) \frac{t^n}{n!}.
\end{split}\end{equation}

Thus, by \eqref{37}, we get
\begin{equation}\begin{split}\label{38}
B_{n,k}^{(\lambda )} (x_1,\cdots,x_{n-k+1}) = B_{n,k}( (1)_{1,\lambda }x_1, (1)_{2,\lambda }x_2, \cdots, (1)_{n-k+1,\lambda }x_{n-k+1}),
\end{split}\end{equation}
where $n \geq k \geq 0$.

Note that
\begin{equation*}\begin{split}
B_{n,k}(1,1,\cdots,1) = S_2(n,k),
\end{split}\end{equation*}
and
\begin{equation*}\begin{split}
B_{n,k}^{(\lambda )} (1,1,\cdots,1) = S_{2,\lambda }(n,k),\,\,(n \geq k \geq 0).
\end{split}\end{equation*}

Kim defined the partially degenerate Bell polynomials which are given by
\begin{equation}\begin{split}\label{39}
e^{x\big((1+\lambda t)^{\frac{1}{\lambda }}-1\big)} = \sum_{n=0}^\infty Bel_{n,\lambda }(x)  \frac{t^n}{n!},\quad (\text{see}\;[8]).
\end{split}\end{equation}

We note that
\begin{equation}\begin{split}\label{40}
&\exp \left( \sum_{j=1}^\infty x_j (1)_{j,\lambda } \frac{t^j}{j!} \right) = \sum_{k=0}^\infty \frac{1}{k!} \left( \sum_{j=1}^\infty x_j (1)_{j,\lambda } \frac{t^j}{j!} \right)^k \\
&= \sum_{k=0}^\infty \sum_{n=k}^\infty B_{n,k}^{(\lambda )} (x_1,x_2, \cdots x_{n-k+1}) \frac{t^n}{n!}\\
&= \sum_{n=0}^\infty \left( \sum_{k=0}^n B_{n,k}^{(\lambda )} (x_1,x_2,\cdots,x_{n-k+1}) \right)       \frac{t^n}{n!}.
\end{split}\end{equation}

Therefore, by \eqref{12} and \eqref{40}, we get
\begin{equation}\begin{split}\label{41}
Bel_n^{(\lambda )} (x_1,x_2,\cdots,x_n) = \sum_{k=0}^n B_{n,k}^{(\lambda )} (x_1,x_2,\cdots,x_{n-k+1}).
\end{split}\end{equation}

Note that
\begin{equation*}\begin{split}
Bel_n^{(\lambda )} (x,x,\cdots,x) &= \sum_{k=0}^n x^k B_{n,k}( (1)_{1,\lambda }, (1)_{2,\lambda },\cdots,(1)_{n-k+1,\lambda } )\\
&= Bel_{n,\lambda }(x),\,\,(n \geq 0).
\end{split}\end{equation*}

Now, we observe that
\begin{equation}\begin{split}\label{42}
&\sum_{n=k}^{\infty} B_{n,k} (0!, 1!\lambda , 2!\lambda ^2, \cdots, (n-k)! \lambda ^{n-k}) \frac{t^n}{n!}\\
&= \frac{1}{k!} \Big( t + \frac{\lambda }{2}t^2 + \frac{\lambda ^2}{3} t^3 + \cdots \Big)^k\\
&= \frac{1}{k!} \Big( - \frac{1}{\lambda}  \log(1-\lambda t) \Big)^k = \frac{1}{k!}  \left( \log \big( \frac{1}{(1-\lambda t)^{\frac{1}{\lambda }}} \big) \right)^k\\
&= \sum_{n=k}^\infty S_{\lambda }(n,k) \frac{t^n}{n!}.
\end{split}\end{equation}

Comparing the coefficients on both sides of \eqref{42}, we have
\begin{equation}\begin{split}\label{43}
B_{n,k} (0!, 1!\lambda , 2!\lambda ^2, \cdots, (n-k)! \lambda ^{n-k}) = S_{\lambda }(n,k),\,\,(n \geq k \geq 0).
\end{split}\end{equation}

\begin{thm}
Let $S_\lambda (n,r)$ be the unsigned degenerate Stirling numbers of the first kind given in \eqref{15} or \eqref{21}. Then, for $n \geq k \geq 0$, we have
\begin{equation}\begin{split}
S_{\lambda }(n,k)=B_{n,k} (0!, 1!\lambda , 2!\lambda ^2, \cdots, (n-k)! \lambda ^{n-k}),
\end{split}\end{equation}
where $B_{n,k} (x_1,x_2,\cdots,x_{n-k+1})$ are the exponential partial Bell polynomials given in \eqref{35}
\end{thm}

\section{Further remarks}

A derangement is a permutation with no fixed points. For example, (2,3,1) and (3,2,1) are derangements of (1,2,3). But (3,2,1) is not because 2 is a fixed point. The number of derangements of an $n$-element set is called the $n$-th derangement number and denoted by $d_n, \;(n \geq 0)$. The derangement number $d_n$ satisfies the following recurrence relation:
\begin{equation}\begin{split}\label{44}
d_n = nd_{n-1} + (-1)^n ,\,\,(n \geq 1).
\end{split}\end{equation}

The generating function of derangement numbers is given by
\begin{equation}\begin{split}\label{45}
\frac{1}{1-t} e^{-t} = \sum_{n=0}^\infty   d_n     \frac{t^n}{n!}.
\end{split}\end{equation}

Thus, by \eqref{45}, we easily get
\begin{equation}\begin{split}\label{46}
d_n = n! \sum_{m=0}^n \frac{(-1)^m}{m!},\,\,(n \geq 0).
\end{split}\end{equation}

For $\lambda \in (0,1)$, the degenerate derangement numbers are defined by the generating function 
\begin{equation}\begin{split}\label{47}
\frac{1}{1-\lambda -t} (1+\lambda t)^{-\frac{1}{\lambda }} = \sum_{n=0}^\infty  d_{n,\lambda }      \frac{t^n}{n!}.
\end{split}\end{equation}

Note that $\lim_{\lambda  \rightarrow 0} d_{n,\lambda } = d_n$, $(n \geq 0)$.
From \eqref{47}, we have
\begin{equation}\begin{split}\label{48}
\sum_{n=0}^\infty  <1>_{n,\lambda }(-1)^n      \frac{t^n}{n!} &= (1+\lambda t)^{-\frac{1}{\lambda }} = \left( \sum_{m=0}^\infty d_{m,\lambda } \frac{t^m}{m!} \right) (1-\lambda -t)\\
&= (1-\lambda ) \sum_{n=0}^\infty   d_{n,\lambda }     \frac{t^n}{n!} - \sum_{n=1}^\infty       n d_{n-1,\lambda } \frac{t^n}{n!}.
\end{split}\end{equation}

By \eqref{48}, we easily get
\begin{equation}\begin{split}\label{49}
d_{0,\lambda } = \frac{1}{1-\lambda },\,\,(1-\lambda )d_{n,\lambda } = (-1)^n <1>_{n,\lambda } + nd_{n-1,\lambda },\,\,(n \geq 1),
\end{split}\end{equation}

and
\begin{equation}\begin{split}\label{50}
(1-\lambda )d_{n+1,\lambda } = (n+\lambda )d_{n,\lambda } + n d_{n-1,\lambda }+ n\lambda (-1)^{n-1} <1>_{n,\lambda } ,\,\,(n \geq 0).
\end{split}\end{equation}

For $\lambda  \in (0,\infty)$ and $\alpha(>0) \in \mathbb{R}$, the degenerate gamma function is defined as
\begin{equation}\begin{split}\label{51}
\Gamma_\lambda  (\alpha) = \int_0^\infty (1+\lambda t)^{-\frac{1}{\lambda }}t^{\alpha-1} dt,\quad (\text{see}\;[7]).
\end{split}\end{equation}

Thus, by \eqref{51}, we get
\begin{equation}\begin{split}\label{52}
\Gamma_\lambda (\alpha+1) = \frac{\alpha}{(1-\lambda )^{\alpha-1}} \Gamma_{\frac{\lambda }{1-\lambda }} (\alpha),\quad (\text{see}\;[7]),
\end{split}\end{equation}
where $\lambda  \in (0,1)$ and $0<\alpha< \frac{1-\lambda }{\lambda }$.

From \eqref{52}, we note that
\begin{equation}\begin{split}\label{53}
\Gamma_\lambda (k) = \frac{(k-1)!}{(1-\lambda )(1-2\lambda )\cdots(1-k\lambda )},\,\,(k \in \mathbb{N},\,\,\lambda  \in (0, \tfrac{1}{k}) ).
\end{split}\end{equation}

Let  $\lambda \in (0,\infty)$. Then $X_{\lambda}$ is the degenerate gamma random variable with  parameters $\alpha (>0)$, $\beta(>0)$, if its probability density function $f$ is given by
\begin{equation}\begin{split}\label{54}
f(x) = \begin{cases}
\frac{1}{\Gamma_\lambda (\alpha)} \beta (\beta x)^{\alpha-1} (1+\lambda x)^{-\frac{1}{\lambda }},&\text{if}\,\, x \geq 0, \\ 
0,&\text{otherwise}.
\end{cases}
\end{split}\end{equation}

Now, we observe that
\begin{equation}\begin{split}\label{55}
(1+\lambda x)^{\frac{t}{\lambda }} &= \sum_{k=0}^\infty \lambda ^{-k} \frac{t^k}{k!} \Big( \log(1+\lambda x)\big)^k \\
&= \sum_{k=0}^\infty t^k	\sum_{n=k}^\infty S_{1,\lambda }(n,k) \frac{x^n}{n!}.
\end{split}\end{equation}

Assume that $X=X_{\lambda}$ is the degenerate gamma random variable with parameters 1,1. Then the $k$-th moment of $X$ is given by
\begin{equation}\begin{split}\label{56}
E[X^k] &= \int_0^\infty x^k \frac{1}{\Gamma_\lambda (1)} (1+\lambda x)^{-\frac{1}{\lambda }} dx\\
&= \frac{1}{\Gamma_\lambda (1)} \Gamma_\lambda  (k+1),\,\,(k \in \mathbb{N}).
\end{split}\end{equation}

For the following discussion, we assume that $\lambda \in (0,\frac{1}{k+1})$, and that $t<1-\lambda$.

On the one hand, the expectation of $(1+\lambda X)^{\frac{t}{\lambda }}$ is given by
\begin{equation}\begin{split}\label{57}
&E[(1+\lambda X)^{\frac{t}{\lambda }}] = \sum_{n=0}^\infty  \left( \sum_{k=n}^\infty    S_{1,\lambda }(k,n) \frac{1}{k!} E[X^k] \right) t^n \\
&= \sum_{n=0}^\infty \left( \sum_{k=n}^\infty S_{1,\lambda }(k,n) \frac{1}{k!} \left( \frac{\Gamma_\lambda (k+1)}{\Gamma_\lambda (1)} \right) \right) t^n\\
&= \sum_{n=0}^\infty \left( \sum_{k=n}^\infty S_{1,\lambda }(k,n) \frac{1}{k!} \times \frac{k!}{\Gamma_\lambda (1) (1-\lambda )(1-2\lambda )\cdots(1-(k+1)\lambda )} \right) t^n\\
&=\sum_{n=0}^\infty \left( \sum_{k=n}^\infty \frac{S_{1,\lambda }(k,n)}{(1-2\lambda )(1-3\lambda )\cdots(1-(k+1)\lambda )} \right) t^n.
\end{split}\end{equation}

\noindent On the other hand,
\begin{equation}\begin{split}\label{58}
&E[(1+\lambda X)^{\frac{t}{\lambda }}] = \int_0^\infty (1+\lambda x)^{\frac{t}{\lambda }} \frac{1}{\Gamma_\lambda (1)} (1+\lambda x)^{-\frac{1}{\lambda }} dx\\
&= \frac{1}{\Gamma_\lambda (1)} \int_0^\infty (1+\lambda x)^{-\frac{1}{\lambda }(1-t)} dx = \frac{1}{\Gamma_\lambda (1)} \frac{1}{1-\lambda -t} \\
&=\frac{1}{\Gamma_\lambda (1)} \frac{1}{1-\lambda -t} (1+\lambda t)^{-\frac{1}{\lambda }} \cdot (1+\lambda t)^{\frac{1}{\lambda }} \\
&= \frac{1}{\Gamma_\lambda (1)} \left( \sum_{l=0}^\infty d_{l,\lambda } \frac{t^l}{l!} \right) \left( \sum_{m=0}^\infty (1)_{m,\lambda } \frac{t^m}{m!} \right)\\
&= \frac{1}{\Gamma_\lambda (1)} \sum_{n=0}^\infty  \left( \sum_{l=0}^n {n \choose l} d_{l,\lambda }(1)_{n-l,\lambda } \right)      \frac{t^n}{n!}.
\end{split}\end{equation}

Therefore, by \eqref{57} and \eqref{58}, we obtain the following theorem.

\begin{thm}
Let  $S_{1,\lambda }(n,k)$ be the degenerate Stirling numbers of the first kind given in \eqref{09} or \eqref{18}, and let $d_{l,\lambda}$ be the degenerate derangement numbers given in \eqref{47}. Then we have
\begin{equation}\begin{split}
 \sum_{k=n}^\infty \frac{S_{1,\lambda }(k,n)n!}{(1-\lambda )(1-2\lambda )\cdots(1-(k+1)\lambda )} = \sum_{l=0}^n {n \choose l} d_{l,\lambda }(1)_{n-l,\lambda },\nonumber
\end{split}\end{equation}
\end{thm}

\noindent where $n \geq 0$, $k \in \mathbb{N}$, and $\lambda  \in (0, \tfrac{1}{	k+1})$.

From \eqref{47}, we have
\begin{equation}\begin{split}\label{60}
\sum_{n=0}^\infty    d_{n,\lambda }    \frac{t^n}{n!} &= \frac{1}{1-\lambda -t} (1+\lambda t)^{-\frac{1}{\lambda }} = \frac{1}{1-\lambda } \left( \frac{1}{1-\frac{t}{1-\lambda }} \right) (1+\lambda t)^{-\frac{1}{\lambda }} \\
&=\frac{1}{1-\lambda } \sum_{m=0}^\infty \left( \frac{t}{1-\lambda } \right)^m \sum_{l=0}^\infty (-1)^l <1>_{l,\lambda } \frac{t^l}{l!}\\
&= \frac{1}{1-\lambda } \sum_{n=0}^\infty \left( n! \sum_{l=0}^n \left( \frac{1}{1-\lambda } \right)^{n-l} (-1)^l <1>_{l,\lambda } \frac{1}{l!} \right) \frac{t^n}{n!}\\
&= \sum_{n=0}^\infty  \left( n! \sum_{l=0}^n \left( \frac{1}{1-\lambda } \right)^{n-l+1}  \frac{(-1)^l}{l!} <1>_{l,\lambda } \right)     \frac{t^n}{n!}.
\end{split}\end{equation} 

Comparing the coefficients on both sides of \eqref{60}, we have
\begin{equation}\begin{split}\label{61}
d_{n,\lambda }   = n! \sum_{l=0}^n \left( \frac{1}{1-\lambda } \right)^{n-l+1}  \frac{(-1)^l}{l!} <1>_{l,\lambda }.
\end{split}\end{equation}

Thus we get the following theorem.

\begin{thm} 
Let $d_{l,\lambda}$ be the degenerate derangement numbers given in \eqref{47}. Then, for $n \geq 0$, we have
\begin{equation}\begin{split}
d_{n,\lambda } = n! \sum_{l=0}^n \left( \frac{1}{1-\lambda } \right)^{n-l+1}  \frac{(-1)^l}{l!} <1>_{l,\lambda }.\nonumber
\end{split}\end{equation}
\end{thm}


\begin{thebibliography}{0}

\bibitem{01}L. Carlitz,
\textit{Degenerate Stirling, Bernoulli and Eulerian numbers,
}Utilitas Math.
{\bf{15}} (1979), 51-88.

\bibitem{02}B. S. El-Desouky, A. Mustafa,
\textit{New results on higher-order Daehee and Bernoulli numbers and polynomials
}Adv. Difference Equ.
{\bf{2016}}, 2016:32, 21 pp.

\bibitem{03}Y.  He, J. Pan, 
\textit{Some recursion formulae for the number of derangements and Bell numbers,
}J. Math. Res. Appl.
{\bf{36}} (2016), no. 1, 15-22.

\bibitem{04}T. Kim,
\textit{Degenerate complete Bell polynomials and numbers,
}Proc. Jangjeon Math. Soc.
{\bf{20}} (2017), no. 4, 533-543.

\bibitem{05}T. Kim,
\textit{$\lambda $-analogue of Stirling numbers of the first kind,
}Adv. Stud. Contemp. Math. (Kyungshang),
{\bf{27}} (2017), no. 3, 423-429.

\bibitem{06}T. Kim, 
\textit{A note on degenerate Stirling polynomials of the second kind,
}Proc. Jangjeon Math. Soc.
{\bf{20}} (2017), no. 3, 319-331.

\bibitem{07}T. Kim, D. S. Kim,
\textit{Degenerate Laplace transform and degenerate gamma function,
}Russ. J. Math. Phys.
{\bf{24}} (2017), no. 2, 241-248.

\bibitem{08}T. Kim, D. S. Kim, G.-W. Jang,
\textit{Extended Stirling polynomials of the second kind and extended Bell polynomials 
}Proc. Jangjeon Math. Soc.
{\bf{20}} (2017), no. 3, 365-376.

\bibitem{09}T. Kim, D. S. Kim, D. V. Dolgy,
\textit{On partially degenerate Bell numbers and polynomials,
}Proc. Jangjeon Math. Soc.
{\bf{20}} (2017), no. 3, 337-345.

\bibitem{10}K. S. K\"{o}lbig,
\textit{The complete Bell polynomials for certain arguments in terms of Stirling numbers of the first kind,
}J. Comput. Appl. Math.
{\bf{51}} (1994), no. 1, 113-116.

\bibitem{11}K. S. K\"{o}lbig, W. Strampp,
\textit{Some infinite integrals with powers of logarithms and the complete Bell polynomials
}J. Comput. Appl. Math.
{\bf{69}} (1996), no. 1, 39-47.

\bibitem{12}J.-W. Park, B. M. Kim, J. Kwon,
\textit{On a modified degenerate Daehee polynomials and numbers,
}J. Nonlinear Sci. Appl
{\bf{10}} (2017), no. 3, 1108-1115.

\bibitem{13}S. Roman,
\textit{The Umbral calculus,
}Pure and Applied Mathematics, 111. Academic Press, Inc. [Harcourt Brace Jovanovich, Publishers], New York, 1984. x+193 pp. ISBN:0-12-594380-6
{\bf{69}} (1996), no. 1, 39-47.

\bibitem{14}Y. Simsek,
\textit{Identities on the Changhee numbers and Apostol-type Daehee polynomials,
}Adv. Stud. Contemp. Math. (Kyungshang),
{\bf{27}} (2017), no. 2, 199-212.
\end{thebibliography}
\end{document}